\documentclass[11pt]{article}
\usepackage[letterpaper,width=135mm,height=203mm]{geometry}

\usepackage{graphicx}
\usepackage{latexsym}

\newtheorem{theorem}{Theorem}
\newtheorem{lemma}[theorem]{Lemma}
\newtheorem{corollary}[theorem]{Corollary}

\newenvironment{Proof}{\textsc{Proof}.}{\hspace{8mm}\QEDmark\smallskip}
\newcommand{\QEDmark}{\mbox{$\Box$}}

\newcommand{\IND}{\alpha}
\newcommand{\TAU}{\tau}

\newcommand{\prism}[1]{{#1\rule{0pt}{8pt}}^{\,p}}
\newcommand{\prismSub}[2]{{#1\rule{0pt}{9pt}_#2^{\,p}}}

\begin{document}


\begin{center}
{\Large \textbf{On the Toughness of  Regular Graphs and Prisms}} \\[2mm]
Geoffrey Boyer \& Wayne Goddard, Clemson University
\end{center}

\begin{abstract}
We contribute results on $r$-regular graphs that do and don't have the maximum possible toughness, namely $r/2$. 
Doty and Ferland showed the existence of a $5$-regular graph with toughness $5/2$ for all even orders except $n= 18$.
Using a computer search we show that there does not exist such a graph for $n=18$.
Also, we provide the first family of $4$-regular
graphs with toughness $2$ that contains claws. For
the prism $G \Box K_2$ of a graph~$G$, we provide several bounds including a sufficient condition for the prism to have the same toughness
as~$G$. In particular, we show that if $G$ has toughness $t\le \frac{1}{2}$ then its prism has toughness $2t$; further,   
the prism of any $r$-regular $r$-connected inflation has toughness~$r/2$ (despite being $(r+1)$-regular) and in general
the prism of any $3$-regular graph has toughness at most~$3/2$. 
\end{abstract}

\section{Introduction}

In 1973, Chv\'{a}tal \cite{chvatal} defined the \textit{toughness} of a graph $G$ to be the minimum value of
$|S|/k(G-S)$ where $k(G-S)$ denotes the number of components of
$G-S$ and the minimum is taken over all cut-sets $S\subseteq V(G)$.
For example, it is immediate that the toughness is at most half the connectivity,
and that every hamiltonian graph has toughness at least $1$. The relationship between toughness
and properties related to hamiltonicity remains an ongoing area of research. In this paper
we instead focus on constructions and bounds on the toughness itself. 
In this regard, a fundamental bound was given by Matthews and Sumner~\cite{MS} 
who showed that if the graph is claw-free, where a \textit{claw} is
an induced copy of $K_{1,3}$, then its toughness equals half its connectivity. 
For a survey of toughness at the time, see~\cite{survey}.

Most of our results are about regular graphs.
An $r$-regular graph is called \textit{supertough} if it has toughness $r/2$.
In the original paper Chv\'{a}tal~\cite{chvatal} noted that 
for $r$ even the power of the cycle $C^{r/2}_n$ (sometimes called a Harary graph) is supertough; in particular,
for $r$ even supertough graphs exist for all orders.
As regards $3$-regular graphs, 
Chv\'{a}tal provided examples (called inflations) of supertough $3$-regular graphs when $n$ is a multiple of $6$
and showed that when $n$ is not a multiple of $6$ such graphs do not exist.
Later Jackson and Katerinis~\cite{JK} showed that being claw-free is also necessary for 
a $3$-regular graph to be supertough. 

A natural place to look for supertough graphs is line graphs, but these can produce only some orders.
In a series of papers, Doty and Ferland~\cite{doty,DFaust05,DFjcmcc08,DFaus09}  provided (inter alia) several constructions of $5$-regular supertough graphs 
showing collectively that there exists a $5$-regular supertough graph
for all even $n$ except for $n=18$. In Section~\ref{s:eighteen} we confirm the nonexistence for $n=18$.
For odd degree $7$ or more, it remains the case that almost nothing is known about orders for which line graphs don't exist.

For even degrees, the powers of cycles noted by Chv\'{a}tal as well as all suitable line graphs all are claw-free. Doty and Ferland~\cite{DFaust05} provided the
first examples of supertough regular graphs with claws. In Section~\ref{s:4reg} we provide the first $4$-regular examples.

Another place to look for supertough graphs is to consider graph products. One natural candidate
is the \textit{prism} $\prism{G}$ of a graph $G$ defined as its cartesian product with $K_2$, or equivalently, taking 
two copies of $G$ and adding a perfect matching between corresponding vertices, since the prism of an $r$-connected
$r$-regular graph is automatically $(r+1)$-connected and $(r+1)$-regular. 
The relationship between toughness and the hamiltonicity of the prism has been of interest recently; see for example~\cite{ESS,SHprism}. 
But we consider here bounds on the toughness of the prism. The first general bounds were provided by 
Casablanca, Diánez, and García-Vázquez~\cite{CDGproceed}.
Earlier, Chao and Han~\cite{CHcycle} established a conjecture of Piazza et al.~\cite{PRSvulnerability}
that a cycle permutation graph (that is, any graph formed by taking two disjoint cycles of the same length and 
adding a perfect matching between them) has toughness at most $4/3$ except when the cycles are triangles. (See also~\cite{GoddardCubic}.)

 In Section~\ref{s:prismLower}
we provide bounds on the toughness of prisms and a characterization on when a graph and its prism have equal toughness. 
Using these ideas, in Section~\ref{s:prismRegular} we  show that many regular prisms are not supertough,
including that the prism of an $r/2$-tough inflation has itself toughness $r/2$.
In Section~\ref{s:double} we consider graphs with small toughness and 
show that if $G$ has toughness $t\le \frac{1}{2}$ then its prism has toughness $2t$.

\section{Supertough $5$-Regular Graphs} \label{s:eighteen}

Despite $5$-regular supertough graphs  existing for all even $n\neq 18$, we show:
 
\begin{theorem}
There is no $5$-regular supertough graph of order $18$.
\end{theorem}

The proof is by exhaustive computer search. We do not know a by-hand proof. We made use of the code \textsf{genreg} provided by Meringer (on github), 
and described in~\cite{Meringer}. We generated all 2,807,105,250,897 connected $5$-regular
graphs of order~$18$, and verified that none of them has the desired toughness. Specifically, we noted that a $5$-regular supertough graph
must have independence number $4$. 

\begin{lemma}
If $G$ is a $5$-regular supertough graph of order $18$ then there is no independent set of size $5$.
\end{lemma}
\begin{Proof}
Suppose $J$ is an independent set of size $5$. There are exactly $25$ edges
from $J$ to $V-J$. Since $|V-J|=13$, there exists some vertex $v\in V-J$
that has at most one neighbor in $J$. It follows that $J \cup \{v \}$ induces
a subgraph with at least five components. By considering $S = V - (J \cup \{v \})$,
it follows that $\tau(G) \le (18-6)/5 < 5/2$, a contradiction.
\end{Proof}

Thus each generated graph was immediately tested for an independent set of size $5$.
There were only 624 graphs with independence number $4$. (It seems likely that this collection could be generated in a more efficient manner.)
The toughness of each of these graphs was then calculated. 

Doty and Ferland~\cite{DFjcmcc08} provided a $5$-regular graph of order $18$ with toughness~$12/5$. It follows that this value is the 
maximum possible toughness.

\section{2-Tough 4-Regular Graphs with Claws}  \label{s:4reg}

Several authors have opined on the existence of supertough graphs.
In the original paper, 
 after noting that the order of a $3$-regular supertough graph is a multiple of $3$ (except $K_4$),
Chv\'{a}tal~\cite{chvatal} expressed the opinion that this behavior was likely for odd $r$ and order sufficiently large.
This was shown to be false by Doty~\cite{doty}.
Noting that the known supertough graphs were claw-free, we wrote about the 
question for $r$-regular graphs for larger $r$:
in \cite{first} we conjectured
that only claw-free graphs could be supertough, while in 
\cite{GoddardCubic} we expressed the opposite belief that almost all $r$-regular graphs are supertough.

It turns out that we were wrong in both cases. Doty and Ferland~\cite{DFaust05} gave the first example of an $r$-regular graph that 
has toughness $r/2$ and claws, and in~\cite{DFjcmcc08} they provided
an infinite family for $r=5$. On the other hand, supertough $r$-regular graphs have independence number at most $2n/(r+2)$, 
but Bollob\'as~\cite{bollobas}  showed that the independence number of a random $r$-regular graph
is at least of the order of $ n \log r / r$, and hence most regular graphs are not supertough. 

In this section we provide a construction of $4$-regular supertough graphs that have claws.
We start with the comment that computer search reveals there are two supertough $4$-regular graphs of order $10$
with claws. Here is one of them. (Note that this is an example of what Alspach and Parsons~\cite{APmeta} called a metacirculant.) 

\begin{center}
\includegraphics[scale=0.7]{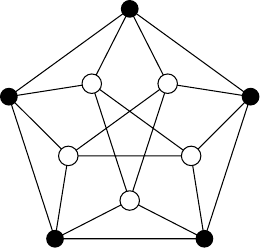}
\end{center}

We produce an infinite family as follows.
For $m\ge 3$, define a graph $J_m$ on $3m-1$ vertices as follows.
Take two disjoint copies of the $m$-cycle, say with vertices $A = \{ a_1, \ldots, a_m \}$ 
and $B = \{ b_1, \ldots, b_m \}$.
Then add vertices $C = \{ c_1, \ldots, c_{m-1} \}$ and join each $c_i$ to each of $a_i,a_{i+1},b_i,b_{i+1}$.
Finally, add edges $a_1 b_1$ and $a_m b_m$. The result is $4$-regular. The graph $J_5$ is shown here.
Note that for $m\ge 4$ the vertices $a_1,a_m,b_1,b_m$ are centers of claws.

\begin{center}
\includegraphics[scale=1.2]{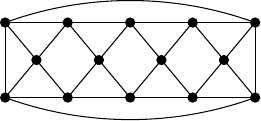}
\end{center}

\begin{lemma}
(a) The graph $J_m$ has connectivity $4$. \\
(b) For $m \ge 5$, a cut-set of size $4$ is either $N(x)$ for some $x \in A \cup B$ or is of the form
$\{ a_i,a_j,b_i,b_j \}$.\\
(c) For $m$ odd, the graph $J_m$ has independence number $m-1$. 
\end{lemma}
\begin{Proof}
(a,b) Consider a cut-set $S$ of size at most $4$. 
Assume first that $S$ does not contain two vertices of $A$. Then the remaining vertices of $A$ in $J_m-S$ lie in one component, say $Z$, and any 
remaining vertex of $C$ is part of $Z$. So there must be a component, say $X$, completely contained in $B$. Every vertex of $A \cup C$ adjacent
to $X$ must be in $S$. It follows readily that $S$ must also contain two vertices of $B$ and that $X$ is a singleton. That is, $S=N(x)$ for some $x\in B$.

The situation where $S$ does not contain two vertices of $B$ is similar. So
suppose that $S$ contains two vertices of both $A$ and $B$. Then it contains no vertex of $C$ 
and it is readily argued that the vertices of $S\cap A$ and $S\cap B$ must align. Thus
$S$ is of the stated form.

(c) The graph $J_m - \{a_1, b_m \}$ has a spanning subgraph consisting of $m-1$ triangles:
$\{c_1,b_1,b_2\}$, $\{c_2,a_2,a_3\}$, and so on. Thus every independent set of size~$m$ must
contain at least one of $\{a_1,b_m\}$. By a symmetric argument, every such 
independent set must contain one of $\{b_1,a_m\}$; but this is a contradiction of the independence.
\end{Proof}

\begin{theorem} \label{t:jm}
For $m\ge 3$ and odd, the graph $J_m$ has toughness $2$.
\end{theorem}
\begin{Proof}
The value $2$ is trivially an upper bound. So we need to show that the graph is $2$-tough.
By the above lemma the graph $J_m$ is $4$-connected. For $m=3$ the graph is claw-free, and so we are done.
So assume $m \ge 5$. Note that each claw in $J_m$ is centered at one of $X = \{ a_1, a_m, b_1, b_m \}$.

Suppose the toughness is less than $2$. That is, there is a cut-set $S$ such that $J_m-S$ has more than $|S|/2$ components.
Out of all such cut-sets, choose one such that $S$ is as large as possible.
Assume the components of $J_m-S$ are $H_1, \ldots, H_k$.
Define $S_i$ as the set of vertices of $S$ that are adjacent to $H_i$ and let $P = \sum_i |S_i|$. That is,
$P$ is the number of pairs $(s,i)$ such that $s\in S$, $ i \in \{1, \ldots, k \}$ 
and vertex $s$ is adjacent to component $H_i$. 
By connectivity, $| S_i | \ge 4$, and so $P \ge 4 k $.
Note that no vertex in $X$ is the center of an induced $K_{1,4}$.
Thus $P \le 2 (|S| - t) + 3t = 2|S| + t$, where $t$ is the number of vertices of $X$
that are in~$S$ and have neighbors in three components of $J_m-S$.
Since $X$ induces a $4$-cycle, and every claw uses two of the edges of that cycle,
it follows that $t\le 2$. Since $|S| \ge 2k - t/2$ but $|S|/k<2$, it follows that 
$t=2$ and further that $P = 4 k$. In particular, $|S_i|=4$ for all $i$.

Suppose there exists a nontrivial component $H_i$.
If $S_i$ is of the form $N(x)$ for some $x \in A \cup B$, then since $J_m - N[x]$ is connected, 
the component $H_i$ must be $J_m-N[x]$. 
Thus in this case there are exactly two components, $H_i$ and $\{x\}$, and the toughness is $2$.
So, by the above lemma we may assume that $H_i$ is a subgraph
that results when the set $\{a_i,a_j,b_i,b_j\}$ is removed. Such a subgraph has
a vertex $c_\ell$ of degree $2$. By adding to $S$ the two neighbors of $c_\ell$ in $H_i$ we
increase the number of components by $1$, and thereby contradict the maximality of $S$.
That is, every component $H_i$ is an isolated vertex. In other words, $S$ is a vertex cover. 
But by the above lemma, the independence number of $J_m$
is only $m-1$, and so $|S|/k \ge 2m/(m-1)>2$, a contradiction.
\end{Proof}

The above construction was found by starting with a line graph and performing a local adjustment;
maybe this idea works in general. It is unclear what happens
 if one insists that \textit{every} vertex is in a claw.

\section{General Lower  Bounds for Prisms} \label{s:prismLower}

In this section we establish some bounds on the toughness of prisms and consider their sharpness.
For the prism $\prism{G}$ 
we refer to the two copies of $G$ as fibers and denote them by $G_1$
and $G_2$. For a set $S$ of vertices in $\prism{G}$ we denote by 
$S_1$ and $S_2$ its restriction to  $G_1$ and $G_2$ respectively.

It is tempting to conjecture that the toughness of the prism is at least the toughness of the original, but that is not true. 
The key point is that one can use as subset of a fiber a set that is not a cut-set. The first result in this regard 
was presented by Casablanca, Di\'{a}nez, and García-V\'{a}zquez~\cite{CDGproceed} without proof. 
We add the proof.

\begin{theorem}  \label{t:generalPrismLower} \cite{CDGproceed}
Let $G$ be a connected graph of order $n$ and independence number $\IND(G)$. Then
\[
   \TAU (\prism{G}) \ge \min \left\{ \TAU(G), \frac{n}{  \IND(G)+1} \right \} .
\]
\end{theorem}
\begin{Proof}
Let $S$ be a cut-set of the prism where $S_i$ is the subset contained in fiber~$G_i$.
Let $k_i$ be the number of components that have \textbf{at least} one vertex in~$G_i$, and note that $k_1,k_2 \le \IND(G)$. 
Let $\prism{k}$ denote
$k( \prism{G} - S )$. There are three cases. 

(i) Both $S_1$ and $S_2$ are cut-sets of $G$. Then 
\begin{equation}  \label{eq:chain}
      \frac{ |S| }{ \prism{k} } \ge \frac{ |S_1| + |S_2| }{ k_1 + k_2 } \ge 
             \min \left \{  \frac{ |S_1| }{ k_1  } , \frac{ |S_2| }{ k_2  }  \right \}  \ge \TAU (G) .
\end{equation}

(ii) Neither $S_1$ nor $S_2$ is a cut-set of $G$. Then $\prism{G}-S$ has two components. Further, 
$|S| \ge n$, since all edges between the two fibers must be severed. In this case $|S| / \prism{k}  \ge n / (\IND(G)+1)$.

(iii)  Assume $S_1$ say is a cut-set of $G$ but $S_2$ is not. 
If there is a component of $\prism{G} - S$ completely contained in $G_2$,
then all edges between the fibers are severed. So $|S| \ge n$ and $\prism{k}  = k_1+1$ and thus $|S| /  \prism{k}  \ge n / (\IND(G)+1)$.
If there is no such component, then $\prism{k}=k_1$. Further, $|S_2| > 0$ since there must be a component completely contained in $G_1$.
Thus $|S| /  \prism{k}  > |S_1|  / k_1 = \TAU ( G ) $
\end{Proof}

At the same time, we observe that one of the terms in the lower bound of Theorem~\ref{t:generalPrismLower} is also an upper bound:

\begin{theorem}  \label{t:generalPrismUpper}
Let $G$ be a connected graph of order $n$ and independence number~$\IND(G)$. Then 
\[
     \TAU (\prism{G}) \le \frac{n}{ \IND(G)+1 } . 
\]
In particular, if $\TAU(G) \ge n/( \IND(G)+1)$, then $ \TAU (\prism{G})  = n/( \IND(G)+1) $.
\end{theorem}
\begin{Proof}
Let $A$ be a maximum independent set of $G$. Then define in the prism the set $S$ to be the union of $A$ from $G_1$ 
and the complement of $A$ from $G_2$. Thus $|S|=n$. In the subgraph $\prism{G}-S$ each remaining vertex of $G_1$ 
is in a component by itself, and there is at least one component containing vertices of $G_2$. 
Thus $k( \prism{G}-S)  \ge \IND(G)+1$.
\end{Proof}

From Theorem~\ref{t:generalPrismLower} one can deduce that the toughness of the prism is at least $\frac{2}{3} \TAU (G)$ for noncomplete $G$.
This follows since the toughness of $G$ is less than $n/\IND(G)$, and thus $\TAU(\prism{G}) / \TAU(G) \ge \IND(G)/(\IND(G)+1)$.
The ratio $\frac{2}{3}$ is asymptotically sharp. For example, consider the graph $K_n-e$. Then its toughness is $(n-2)/2$, while 
its prism has toughness $n/3$. 

One can also ask when does graph $G$ and its prism have the same toughness. It is immediate that if $G$ is disconnected
then both have toughness $0$. And by Theorem~\ref{t:generalPrismUpper}, equality holds if $\TAU(G) \ge n/( \IND+1)$.
We say a graph $G$ has \textit{complementary tough sets} if there exist tough sets $S_1$ and $S_2$ such that $S_1 \cup S_2 = V(G)$.
For example, a bipartite graph with toughness $1$ has complementary tough sets, namely the two partite sets.

\begin{theorem} \label{t:complementaryTough}
For connected graph $G$ with $\TAU(G) < n/( \IND+1)$,  it holds that $\TAU ( \prism{G} ) = \TAU ( G ) $ if and only if
$G$ has complementary tough sets.
\end{theorem}
\begin{Proof}
Assume $G$ has complementary tough sets $S_1$ and $S_2$.
We form a cut-set of $\prism{G}$ by taking the vertices of $S_1$ in $G_1$ and the vertices of $S_2$ in $G_2$. 
Since $S_1 \cup S_2 = V(G)$, each component of $\prism{G}-S$ is completely contained within one of the fibers
and indeed $G-S$ has $k(G_1-S_1) + k(G_2-S_2)$ components. It follows that $\TAU( \prism{G} ) \le (|S_1|+|S_2|) / (
k(G_1-S_1) + k(G_2-S_2) ) = \TAU (G )$. By Theorem~\ref{t:generalPrismLower} we have equality.

Assume  $\TAU ( \prism{G} ) = \TAU ( G ) $. Let $S$ be a tough set of $\prism{G}$. 
Since $\TAU(\prism{G}) < n /(\IND(G)+1)$, by the proof of Theorem~\ref{t:generalPrismLower} it must be that the restrictions
 $S_1$ and $S_2$ of $S$ 
are both cut-sets in $G_1$ and $G_2$. Thus from the inequality chain~(\ref{eq:chain}), 
it is necessary that  $|S_1| / k( G_1-S_1 ) = |S_2|/k(G_2-S_2) = \TAU(G)$. 
Further, we need $k(\prism{G}-S) = k( G_1-S_1) + k(G_2-S_2)$. Thus 
each component of $k(\prism{G}-S)$ is completely contained within a fiber, or in other words,
every edge of $\prism{G}$ joining the fibers has an end in $S$. It follows that $S_1$ and $S_2$ are complementary tough sets of $G$.
\end{Proof}

As an application of this result, we note for example that it applies to the graph $J_m$ described earlier:

\begin{theorem}  \label{t:prismJm}
For odd $m$ the graph $J_m$ has complementary tough sets and so $\TAU( \prismSub{J}{m} )=2$. 
\end{theorem}
\begin{Proof}
Since Theorem~\ref{t:jm} showed that $\tau(J_m)=2$, 
we need two sets $S_1$ and $S_2$ whose union is the entire vertex set and $|S_1|/k(J_m-S_1) = |S_2|/k(J_m-S_2) = 2$. 
Recall that the vertices of one cycle of $J_m$ are $A = \{ a_1, \ldots, a_m \}$, 
the other $B= \{ b_1, \ldots, b_m \}$, and the remaining vertices $C = \{ c_1, \ldots, c_{m-1} \}$ such that vertex $c_i$ is adjacent to vertices $a_i$, $a_{i+1}$, $b_{i}$, and $b_{i+1}$. 
Define
\[
S_1= A \cup B - \{a_1, b_m \} , \qquad 
S_2=\{a_1, a_3, \ldots, a_{m-2} \} \cup \{b_3, b_5, \ldots, b_m \} \cup C .
\] 
It is immediate that $|S_1|=|S_2|=2(m-1)$ and their union is the whole vertex set. Each vertex of $C$ is in a separate
component of $J_m - S_1$ for a total of $m-1$ components. 
In the subgraph $J_m -S_2$, each odd cycle is broken into $(m-1)/2$ pieces for a total of $m-1$ components.    
\end{Proof}

%
%


\section{Prisms of Regular Graphs} \label{s:prismRegular}

The inflation of a graph $G$ can be defined as $L(S(G))$, that is, the line graph of the subdivision graph.
As noted above, Chv\'{a}tal used the inflation of a graph to produce graphs that are supertough. So one might hope that the 
prism of an inflation is also supertough. However, despite the fact that an inflation has small tough sets (for example the neighborhood of any vertex),
the following holds:

\begin{theorem}
For $r\ge 2$ an $r$-connected $r$-regular inflation $G$ has complementary tough sets.
\end{theorem}
\begin{Proof}
Assume $G=L(S(H))$ where graph $H$ has order $m$. Note that $G$ has order $n=rm$. 
Since $G$ is $2$-connected, it follows
that the graph $S(H)$ and hence the graph $H$ is $2$-edge-connected. 
By Robbins' Theorem
one can orient the edges of $H$ such that the result is strongly connected.
In particular, every vertex is the head and the tail of at least one arc.

Now for each oriented edge $e$ in $H$ take the edge in $S(H)$ corresponding
to the head, and let $X$ denote the resultant set of edges. 
Note that every vertex of~$H$ is in a different component of $S(H)-X$ 
and each of these components is nontrivial. It follows that $k(G-X)=m$ while $|X|=n/2$.
Thence $|X|/k(G-X) = (n/2)/(n/r) = r/2$ and
$X$ is a tough set for $G$. Similarly, the complement of $X$ 
formed from taking each tail is a tough set. Thus $G$ has complementary
tough sets.
\end{Proof}

It follows from Theorem~\ref{t:complementaryTough}  that if $G$ is an inflation of degree $r$ then $\TAU ( \prism{G} ) = r/2$.
We observe next that the upper bound generalizes to all $G$ that are cubic graphs:

\begin{theorem}  \label{t:prismCubic}
If graph $G\neq K_4$ has maximum degree at most $3$, then its prism has toughness at most $3/2$.
\end{theorem}
\begin{Proof}
Consider a proper $3$-coloring $(A_1,A_2,A_3)$ of $G$.
Let $X_1$ be the subset of~$A_3$ whose vertices have at most one neighbor in $A_1$.
Then define set $S$ of the prism as: all vertices of $A_1 \cup X_1$ from the first fiber, and
all vertices of $A_2 \cup (A_3 - X_1)$ from the second fiber. Note that $|S|=n$ where $n$ is the 
order of~$G$.

Consider the subgraph $\prism{G} - S$.
Note that each vertex of $A_2$ in the first fiber $G_1$ is in a separate component, since each vertex
of $A_3 - X_1$ has at most one neighbor in $A_2$. Similarly with each vertex of $A_1$ in the 
second fiber $G_2$. Hence $k( \prism{G} - S ) = |A_1|+|A_2|$. We can choose
this quantity to be at least~$\frac{2}{3} n$. It follows that $\TAU( \prism{G} ) \le n/(2n/3) = 3/2$. 
\end{Proof}

For maximum degree $4$ the proof of Theorem~\ref{t:prismCubic} does not immediately generalize. 
If a vertex has at most one neighbor in each of two of the colors, then it is good;
but it is unclear how to handle vertices that have two neighbors of each of two colors. 
It is unclear whether the toughness of prism of  $4$-regular noncomplete graph can be more than $2$. 
 For some natural candidates, such as the supertough $4$-regular graphs $J_m$ described earlier,  their prism has toughness $2$ (Theorem~\ref{t:prismJm}).
 
 However,
Theorem~\ref{t:prismCubic} does not generalize to $5$-regular graphs. For example, 
there are two supertough $5$-regular graphs of order $12$, namely the icosahedron and the line graph of $K_{3,4}$, and 
it can be shown by computer that 
both of these have a prism whose toughness is~$8/3$. But it is an open question 
whether the prism of a $5$-regular graph can be supertough.


\section{Upper Bounds and Prisms of Graphs with Small Toughness}  \label{s:double}

Casablanca, Di\'{a}nez, and García-V\'{a}zquez in~\cite{CDGproceed} provided an upper bound on the toughness of a prism.
We add the proof.

\begin{theorem} \cite{CDGproceed} \label{t:prismDoubler}
Let $G$ be a connected graph of order $n$ and minimum degree~$\delta(G)$. 
Then 
\[
     \TAU (\prism{G})  \le \min \left \{ 2 \TAU(G) , \frac{\delta(G)+1}{2} \right \} .
 \]
\end{theorem}
\begin{Proof}
To show that $2\TAU(G)$ is an upper bound, take the tough set for $G$ and use it in both fibers. 
To show that $(\delta(G)+1)/2$ is an upper bound, take $S$ to be all neighbors of one vertex of minimum degree.
\end{Proof}


There are graphs where $\TAU( \prism{G} ) = 2 \TAU ( G ) $. 
Consider for example a path: this has toughness $\frac{1}{2}$,
but its prism is hamiltonian and so has toughness $1$. 
Further, a consequence of the main result of Ellingham et al.~\cite{ESS} is that the prism of
every $P_4$-free graph with toughness $\frac{1}{2}$ has toughness $1$.
In fact, this doubling is true of all graphs with toughness at most $\frac{1}{2}$, as we will show. 

For a vertex $v$ in $\prism{G}$, we denote its neighbor in the other $G$-fiber by $v'$. 
Given a tough set $S$ of $\prism{G}$, we say vertex $v \in S$ is \textit{mirrored} if $v' \in S$; 
otherwise $v$ is \textit{split}. We say $S$ is mirrored if every vertex in $S$ is mirrored.

\begin{theorem}
If $\TAU( \prism{G} ) < 1$,  then every tough set of $\prism{G}$ is mirrored. 
\end{theorem}
\begin{Proof} 
We may assume $G$ is connected.  Suppose there exists a tough set $S$ of $\prism{G}$ that is not mirrored.

We mark two things in stages: (a) split vertices, and (b) components of $\prism{G}-S$. At each stage we ensure the following:    
    \begin{itemize}
        \item[(i)] The number of marked split vertices is at least the number of marked components;
        \item[(ii)] For every marked split vertex $w$, the component containing $w'$ is marked; and
        \item[(iii)] The subgraph of $\prism{G}$ induced by the marked vertices and components is connected.
    \end{itemize}
    
 Start by marking some split vertex $v$ and marking all components adjacent to $v$. For each marked component $C$, 
 except the one containing $v'$, 
 and for every vertex $w\in C$ adjacent to $v$, it must be that $w'$ is in $S$ (else $w$ is in the same component as $v'$); 
 mark all such $w'$ for each $C$. The above properties hold. 

 If there is a marked split vertex $x$ that is adjacent to an unmarked component, then proceed as follows. 
 Mark all unmarked components adjacent to $x$. Recall that necessarily the component containing $x'$ is already marked. 
 As above, for each newly marked component $C$ and $w\in C$ adjacent to $x$, it must hold that $w'$ is in~$S$. 
 By the above property it must be that $w'$ is not marked. So mark each such~$w'$ for each $C$. The above properties hold. 
 Repeat until there does not exist a marked split vertex adjacent to an unmarked component. 
    
        If all components of $\prism{G}-S$ are marked, then we have $|S|\geq k( \prism{G} -S)$ and so $\TAU( \prism{G} ) \geq 1$, a contradiction.  
 So assume not all components of $\prism{G}-S$ are marked. From this, it follows that there must be vertices of $S$ that are not marked;
let $S'$ be the set of all such vertices. Note that the subgraph $H$ of $\prism{G}$ induced by the marked vertices and components is connected. 
 Since $H$ is maximal with respect to components of $\prism{G}-S$ adjacent to it, it must be that $S'$ is a cut-set of $\prism{G}$ that leaves $H$ as a component. 
Thus  $k( \prism{G} - S' ) > k( \prism{G} - S) - |S-S'| $. Since $|S| < k( \prism{G} - S) $, 
it follows that  $ |S'| / k( \prism{G} - S' ) < |S|/ k( \prism{G} - S)$. 
This is a contradiction of $S$ being a tough set.
\end{Proof}

\begin{theorem}\label{t:double}
    If $G$ is a graph with toughness $t$, then the toughness of $\prism{G}$ is at least $\min \{ 2t,1 \}$.
\end{theorem}
\begin{Proof}
Assume $\prism{G}$ has toughness less than $1$. Then by the above theorem a tough set $S$ is mirrored.
Let $\bar{S}$ denote the projection of $S$ onto $G$. It follows that $|S| / k( \prism{G} - S ) = 2 |\bar{S}| / k(  G-\bar{S} ) \ge 2 t$.
 \end{Proof}

The following results follow immediately from Theorem~\ref{t:double}. 

\begin{corollary}\label{c:smallDub}
(a)   If $G$ is a graph with toughness $t\leq1/2$, then $\TAU( \prism{G} )=2t$. \\
(b)   If $G$ is bipartite with toughness $t$, then $\TAU( \prism{G} )=\min \{ 2t, 1 \}$.
\end{corollary}
\begin{Proof}
For any graph, $2t$ is an upper bound by Theorem~\ref{t:prismDoubler}.
 Further, the prism of a bipartite graph is bipartite and so has toughness at most $1$.
\end{Proof}




\section{Acknowledgements}

The result in Section 2 was made using the Palmetto 2 cluster, and is based on work supported by the National Science Foundation under Grant Nos. MRI 1228312, II NEW 1405767, MRI 1725573, and MRI 2018069. Any opinions, findings, and conclusions or recommendations expressed in this material are those of the author(s) and do not necessarily reflect the views of the National Science Foundation.

\end{document}